УДК 519.856

# Безградиентные двухточечные методы решения задач стохастической негладкой выпуклой оптимизации при наличии малых шумов не случайной природы

*Баяндина Анастасия Сергевна[1]* [anast.bayandina@gmail.com](anast.bayandina@gmail.com)

*Гасников Александр Владимирович[1,2]* [gasnikov.av@mipt.ru](gasnikov.av@mipt.ru)

*Гулиев Фариман Шакирович[1]* [guliev@phystech.edu](guliev@phystech.edu)

*Лагуновская Анастасия Александровна[1]* [a.lagunovskaya@phystech.edu](a.lagunovskaya@phystech.edu)

[1] Факультет управления и прикладной математики Национального исследовательского Университета «Московский физико-технический институт».
141700, Россия, Московская область, г. Долгопрудный, Институтский переулок, д. 9
[2] Институт проблем передачи информации им. А.А. Харкевича Российской академии наук
127051, Россия, г. Москва, Большой Каретный переулок, д.19 стр. 1

**Аннотация**
В данной работе изучаются негладкие выпуклые задачи стохастической оптимизации с двухточечным оракулом нулевого порядка, т.е. на каждой итерации наблюдению доступны значения реализации функции в двух выбранных точках. Эти задачи предварительно сглаживаются с помощью известной техники двойного сглаживания (Б.Т. Поляк), а затем решаются с помощью стохастического метода зеркального спуска. В данной работе, по-видимому, впервые получены условие на допустимый уровень шума неслучайной природы, проявляющегося при вычислении реализации функции, при котором сохраняется оценка скорости сходимости метода.
**Ключевые слова:** Метод зеркального спуска, шумы, стохастическая оптимизация, безградиентные методы, техника двойного сглаживания.

## 1. Введение

В данной статье изучается специальная техника решения негладких задач выпуклой оптимизации двухточечными безградиентными методами, получившая название "техника двойного сглаживания". Эта техника, восходящая к работам Бориса Теодоровича Поляка (см., например, [1, 2]), до настоящего момента использовалась в основном только в условиях отсутствия шумов неслучайной природы [3, 4]. Однако во многих приложениях важным является наличие таких шумов [2, 5–7]. В настоящей работе будет показано, что негладкость задачи фактически никак не влияет на оценки скорости сходимости оптимальных методов, которые нуждаются в небольшой корректировке (дополнительной рандомизации) в негладком случае, но зато влияет на требование малости шума. В статье дается количественная оценка того, насколько должны быть выше требования к точности вычисления значения функции в негладком случае по сравнению с гладким.

В п. 2 приводится, в немного сказочной форме (для большей наглядности), возможная содержательная постановка задачи, приводящая именно к описанному выше контексту. Формулируется результат о скорости сходимости оптимальных процедур.

В п. 3 излагается вариант метода зеркального спуска с неточным оракулом в специальной форме [7], которая понадобится в дальнейшем (в п. 5) для описания оптимального метода.

В п. 4 описывается техника (однократного) сглаживания негладкой задачи [8]. С помощью этой техники, примененной два раза, в следующем пункте строится сведение задачи с оракулом нулевого порядка, выдающим значения реализации функции в двух точках, к оракулу первого порядка, выдающему стохастический градиент.

В п. 5 с помощью техники двойного сглаживания и метода из п. 3 строится оптимальный метод, работающий с точностью до логарифмических множителей по известным нижним оракульным оценкам [3, 9]. Изучаются требования к шуму не случайной природы, при котором оценки скорости сходимости сохраняют свой вид (с точностью до числовых мультипликативных констант).

## 2. Содержательная интерпретация задачи

Представим себе, что на некоторой воображаемой планете живут "оптимизаторы". Каждый из них стремиться жить с наименьшими потерями нервов. Оптимизаторы очень чувствительные: переживают / нервничают, когда у них что-то не получается. Планета характеризуется выпуклой функцией $f(x)$, аргумента $x \in Q \subseteq \mathbb{R}^n$, где $Q$ – выпуклое множество (множество допустимых стратегий оптимизаторов). Функция $f(x)$ обозначает нервные потери оптимизатора за день, если оптимизатор в этот день использовал стратегию $x \in Q$. Однако значение функции $f(x)$ не наблюдаемы. Оптимизатор, использовавший в $k$-й день стратегию $x^k$, может наблюдать лишь зашумленную реализацию значения этой функции

$$\tilde{f}(x^k, \xi^k) = f(x^k, \xi^k) + \delta(x^k, \xi^k),$$

$$E_{\xi^k}\left[f(x^k, \xi^k)\right] = f(x^k), \; \left|\delta(x^k, \xi^k)\right| \le \delta.$$

Предполагается, что $\{\xi^k\}_k$ – независимые одинаково распределенные (как $\xi$) случайные величины (их реализации одинаковы для всех оптимизаторов), функция $f(x,\xi)$, как функция $x$, является в $\varepsilon/M$-окрестности множества $Q$ ($Q_{\varepsilon/M}$):

- выпуклой, но необязательно гладкой;
- удовлетворяющей условию

$$|f(y,\xi) - f(x,\xi)| \le M\|y - x\|_2. \qquad (1)$$

Успешность оптимизатора, прожившего $N \gg 1$, дней измеряется величиной[1]

$$\text{Regret}_f\left(\{x^k\}_{k=0}^{N-1}\right) = E\left[\frac{1}{N}\sum_{k=0}^{N-1} f(x^k)\right] - \min_{x \in Q} f(x). \qquad (2)$$

Если бы оптимизатор сразу знал "как правильно жить", то оптимальной (с точки зрения минимизации введенного регрета) стратегией была бы стратегия $x^k \equiv x_*$, которая давала бы нулевой регрет (2). Однако в начальный момент оптимизаторы ничего не знают. Все оптимизаторы стартуют с состояния $x^0$. Введем прокс-функцию $d(x) \ge 0$ ($d(x^0) = 0$), которая предполагается сильно выпуклой относительно выбранной нормы $l_p = \|\ \|_p$, $1 \le p \le 2$, с константой сильной выпуклости $\ge 1$ [10]. Положим $R^2 = V(x_*, x^0)$, где прокс-расстояние (расстояние Брэгмана) определяется формулой

$$V(x,z) = d(x) - d(z) - \langle \nabla d(z), x - z\rangle.$$

Если оптимальная стратегия $x_*$ – не единственная, то в определении $R^2$ выбирается такая стратегия $x_*$, которая доставляет минимум $V(x_*, x^0)$.

Рассматривается семейная пара из двух оптимизаторов, которая каждый день может наблюдать $\tilde{f}(x_m^k, \xi^k)$, $\tilde{f}(x_w^k, \xi^k)$, обмениваясь между собой информацией. Таким образом, стратегией семейной пары является способ выбора

$$x_m^k\left(\tilde{f}(x_m^{k-1}, \xi^{k-1}), \tilde{f}(x_w^{k-1}, \xi^{k-1}); ...; \tilde{f}(x_m^0, \xi^0), \tilde{f}(x_w^0, \xi^0)\right),$$

$$x_w^k\left(\tilde{f}(x_m^{k-1}, \xi^{k-1}), \tilde{f}(x_w^{k-1}, \xi^{k-1}); ...; \tilde{f}(x_m^0, \xi^0), \tilde{f}(x_w^0, \xi^0)\right).$$

В данной статье решается задача поиска оптимальной стратегии семейной пары, доставляющей минимально возможное (с точностью до логарифмических множителей) значение регрету каждого члена семьи [6, 7]

---

[1] Обычно эту величину называют псевдо-регретом [11], однако, для краткости далее в этой статье будем использовать название "регрет".

$$\text{Regret}_f\left(\{x^k\}_{k=0}^{N-1}\right) = \tilde{O}\left(n^{1/q}\frac{MR}{\sqrt{N}}\right), \qquad (3)$$

где $1/p + 1/q = 1$, $M$ определяется согласно (1), а запись $A(N,n) = \tilde{O}(B(N,n))$ означает, что существует такая константа $C_q$, возрастающая при $q \in [2, \infty]$ от $O(1)$ до $O(\sqrt{\ln n})$ (см. Лемму 4 в п. 5), что для всех $n$ и $N$ справедливо неравенство $A(N,n) \leq C_q B(N,n)$.

Если дополнительно известно, что функция $f(x)$ является сильно выпуклой в норме $l_2$ с достаточно большой константой $\gamma > 0$, то оценку (3) можно улучшить

$$\text{Regret}_f\left(\{x^k\}_{k=0}^{N-1}\right) = O\left(n\frac{M^2 \ln N}{\gamma N}\right). \qquad (4)$$

Отличие рассматриваемой в данной статье постановки от классической задачи получения оптимальных оценок для выпуклых (вообще говоря, негладких) задач стохастической оптимизации с двухточечным оракулом [3] заключается в наличие небольшого шума $\delta > 0$ неслучайной природы. Это допущение является важным шагом к практической адаптации существующих подходов [1, 5–7]. Далее в статье приводится условие на уровень шума $\delta$, при котором оценки (3), (4) получаются такими же (с точностью до абсолютной мультипликативной константы), как если бы этого шума не было, т.е. имела место несмещенность $\delta = 0$.

Новизна по сравнению с работами [6, 7] заключается в том, что в данной статье для получения оценок (3), (4) не предполагается, что функция $f(x)$ имеет равномерно ограниченную константу Липшица стохастического градиента.[2]

## 3. Метод зеркального спуска для задач стохастической оптимизации

Предположим, что необходимо решать задачу стохастической выпуклой оптимизации [12]

$$F(x) = E_{\tilde{\eta}}\left[F(x, \tilde{\eta})\right] \to \min_{x \in Q}. \qquad (5)$$

На каждой итерации можно один раз обратиться к оракулу за зашумленным значением стохастического градиента $\nabla_x \tilde{F}(x^k, \eta^k)$.

Пусть для любых $\tilde{N} \leq N$ [7] ($\Theta^{k-1}$ – сигма алгебра, порожденная $\eta^1, \ldots, \eta^{k-1}$ [12])

---

[2] Здесь и далее предполагается знакомство читателей с понятием стохастического (суб-)градиента, например, в объеме книги [12].

$$\sup_{\left\{x^k=x^k\left(\xi^1,\ldots,\xi^{k-1}\right)\right\}_{k=1}^{\tilde{N}} \in Q_{\varepsilon/M}} E\left[\frac{1}{\tilde{N}}\sum_{k=1}^{\tilde{N}}\left\langle E_{\eta^k}\left[\nabla_x F\left(x^k,\eta^k\right)-\nabla_x \tilde{F}\left(x^k,\eta^k\right)\Big|\Theta^{k-1}\right], x^k - x_*\right\rangle\right] \leq \sigma, \qquad (6)$$

где $x_*$ – решение задачи (5) (если решение $x_*$ – не единственно, то выбираем такое решение $x_*$, которое доставляет минимум $V(x_*, x^0)$),

$$E_{\eta^k}\left[\nabla_x F\left(x^k,\eta^k\right)\right] = \nabla F\left(x^k\right),\ E_{\eta^k}\left[\left\|\nabla_x \tilde{F}\left(x^k,\eta^k\right)\right\|_q^2\right] \leq \tilde{M}^2, \qquad (7)$$

$\left\{\eta^k\right\}_{k=0}^{N-1}$ – независимые одинаково распределенные случайные величины.

Опишем шаг метода зеркального спуска [7, 10] (МЗС)

$$x^{k+1} = \operatorname{Mirr}_{x^k}\left(h\nabla_x \tilde{F}\left(x^k,\eta^k\right)\right),\ \operatorname{Mirr}_{x^k}(\mathbf{v}) = \arg\min_{x \in Q}\left\{\left\langle \mathbf{v}, x-x^k\right\rangle + V\left(x,x^k\right)\right\}, \qquad (8)$$

размер шага $h$ будет выбран позже (см. формулу (11)).

Основное свойство метода [7, 10]
$$2V\left(x,x^{k+1}\right) \leq 2V\left(x,x^k\right) + 2h\left\langle \nabla_x \tilde{F}\left(x^k,\eta^k\right), x-x^k\right\rangle + h^2\left\|\nabla_x \tilde{F}\left(x^k,\eta^k\right)\right\|_q^2. \qquad (9)$$

Поскольку $\left\{\eta^k\right\}_{k=0}^{N-1}$ – независимые одинаково распределенные случайные величины, то из (9) можно получить

$$F\left(x^k\right) - F(x) \leq \left\langle \nabla F\left(x^k\right), x^k - x\right\rangle = \left\langle E_{\eta^k}\left[\nabla_x F\left(x^k,\eta^k\right)\Big|\Theta^{k-1}\right], x^k - x\right\rangle \leq$$
$$\leq \left\langle E_{\eta^k}\left[\nabla_x F\left(x^k,\eta^k\right)\Big|\Theta^{k-1}\right] - \nabla_x \tilde{F}\left(x^k,\eta^k\right), x^k - x\right\rangle +$$
$$+ \frac{1}{h}\left(V\left(x,x^k\right) - V\left(x,x^{k+1}\right)\right) + \frac{h}{2}\left\|\nabla_x \tilde{F}\left(x^k,\eta^k\right)\right\|_q^2. \qquad (10)$$

Беря условное математическое ожидание $E_{\eta^k}\left[\ \cdot\ \Big|\Theta^{k-1}\right]$ от обеих частей неравенства (10), получим

$$F\left(x^k\right) - F(x) \leq \left\langle \nabla F\left(x^k\right), x^k - x\right\rangle \leq$$
$$\leq \left\langle E_{\eta^k}\left[\nabla_x F\left(x^k,\eta^k\right) - \nabla_x \tilde{F}\left(x^k,\eta^k\right)\Big|\Theta^{k-1}\right], x^k - x\right\rangle +$$
$$+ \frac{1}{h}\left(V\left(x,x^k\right) - E\left[V\left(x,x^{k+1}\right)\Big|\Theta^{k-1}\right]\right) + \frac{h}{2}\underbrace{E_{\eta^k}\left[\left\|\nabla_x \tilde{F}\left(x^k,\eta^k\right)\right\|_*^2\Big|\Theta^{k-1}\right]}_{\leq \tilde{M}^2}.$$

Если просуммировать последнее неравенство по $k = 0,\ldots,N-1$, а затем взять полное математическое ожидание от обеих частей, положив $x = x_*$ (если решение задачи (5) $x_*$ – не единственно, то выбираем такое решение $x_*$, которое доставляет минимум $V(x_*, x^0)$), то получим

$$\text{Regret}_F\left(\left\{x^k\right\}_{k=0}^{N-1}\right) \le (hN)^{-1} V(x_*, x^0) + \tilde{M}^2 h/2 + \sigma \le \sqrt{\frac{2\tilde{M}^2 R^2}{N}} + \sigma,$$

где

$$R^2 = V(x_*, x^0), \quad h = \frac{R}{\tilde{M}}\sqrt{\frac{2}{N}} = \frac{\varepsilon}{\tilde{M}^2}, \qquad (11)$$

т.е.

$$\text{Regret}_F\left(\left\{x^k\right\}_{k=0}^{N-1}\right) \le \tilde{M} R \sqrt{\frac{2}{N}} + \sigma. \qquad (12)$$

Другими словами, после

$$N(\varepsilon) = \frac{2\tilde{M}^2 R^2}{\varepsilon^2} \qquad (13)$$

дней (итераций)

$$\text{Regret}_F\left(\left\{x^k\right\}_{k=0}^{N-1}\right) \le \varepsilon + \sigma. \qquad (14)$$

Оценки (12)–(14) оптимальны с точностью до небольших мультипликативной числовых констант [8]. Используя вместо МЗС метод двойственных усреднений (МДУ) [13], можно прийти к аналогичным оценкам, но уже с адаптивным выбором шагов (11), т.е. $h$ можно выбирать не зависящим от числа итераций $N$ (желаемой точности $\varepsilon$) образом. В случае, когда оценка $\tilde{M}$ не известна, можно использовать вариант МЗС из обзора [14] – это распространяется и на МДУ. Можно перенести метод и на композитные постановки задач [10]. Все написанное выше с сохранением оценок и метода также переносится на задачи онлайн оптимизации [7, 11].

Проведенные рассуждения позволяют (подобно [13]) попутно получить, что

$$E\left[\frac{1}{2}\|x_* - x^k\|^2\right] \le E\left[V(x_*, x^k)\right] \le 2V(x_*, x^0) = 2R^2, \quad k = 0,\dots,N. \qquad (15)$$

Если дополнительно известно, что функция $F(x)$ является сильно выпуклой в норме $l_2$ с достаточно большой константой $\gamma > 0$, то МЗС (8) можно адаптировать под эту специфику

$$x^{k+1} = \text{Mirr}_{x^k}\left(h_{k+1} \nabla_x \tilde{F}(x^k, \eta^k)\right), \quad h_k = (\gamma \cdot k)^{-1}. \qquad (16)$$

Тогда оценка (12) перепишется следующим образом ($\tilde{M}$ посчитано в $l_2$-норме)

$$\text{Regret}_F\left(\left\{x^k\right\}_{k=0}^{N-1}\right) \le \frac{\tilde{M}^2}{2\gamma N}(1 + \ln N) + \sigma. \qquad (17)$$

Оценка (17) оптимальна с точностью до $\sim \ln N$ [8]. Впрочем, нам не известно устраним ли этот логарифмический множитель. Для задач обычной стохастической сильно выпуклой оптимизации, в которых критерием качества является не регрет, а просто итоговая невязка по функции $E\left[F(x^N)\right] - \min_{x \in Q} F(x)$ известно, что этот множитель может быть устранен, например, с помощью процедуры рестартов [15]. Для более общих задач стохастической сильно выпуклой онлайн оптимизации, включающих рассматриваемый класс задач с регретом в качестве критерия качества, этот логарифмический множитель устранить уже нельзя [16].

Можно показать, что все последующие рассуждения можно проводить, беря за основу одно из отмеченных здесь обобщений МЗС. Причем все оценки можно писать в категориях вероятностей больших уклонений [12, 17], а не в среднем, как сейчас. Однако в этой статье мы ограничимся простейшим вариантом, чтобы лучше пояснить основную схему перенесения результатов пп. 3, 4 на случай, когда доступны только зашумленные реализации оптимизируемой функции п. 2.

### 4. Сглаживание задачи

Изложим схему сглаживания, позволяющую свести (см. п. 5) постановку п. 2 к постановке п. 3. Тогда можно будет воспользоваться формулами (5)–(8), (11)–(17).

Пусть $\tilde{e} \in RB_2^n(1)$ – случайный вектор, равномерно распределенный на шаре единичного радиуса в норме $l_2$ в $\mathbb{R}^n$. Сгладим (следуя [1, 3, 5–8]) исходную функцию $F(x,\tilde{\eta})$ с помощью локального усреднения по шару радиуса $\mu > 0$ (см. (19)), которое будет выбрано позже ($\tilde{\eta}$ и $\tilde{e}$ предполагаются независимыми случайными величинами),

$$F^\mu(x,\tilde{\eta}) = E_{\tilde{e}}\left[F(x+\mu\tilde{e},\tilde{\eta})\right], \; F^\mu(x) = E_{\tilde{e},\tilde{\eta}}\left[F(x+\mu\tilde{e},\tilde{\eta})\right].$$

Заменим исходную задачу (5) следующей задачей

$$F^\mu(x) = E_{\tilde{\eta}}\left[F^\mu(x,\tilde{\eta})\right] \to \min_{x \in Q}. \tag{18}$$

Легко проверить, что если (см., например, (1))

$$\left|F(y,\tilde{\eta}) - F(x,\tilde{\eta})\right| \le M\|y-x\|_2,$$

то

$$0 \le F^\mu(x,\tilde{\eta}) - F(x,\tilde{\eta}) \le M\mu, \; 0 \le F^\mu(x) - F(x) \le M\mu.$$

Поэтому, если

$$M\mu \le \varepsilon/2, \tag{19}$$

и

$$\operatorname{Regret}_{F^\mu}\left(\{x^k\}_{k=0}^{N-1}\right) \le \varepsilon/2,$$

то

$$\operatorname{Regret}_F\left(\{x^k\}_{k=0}^{N-1}\right) = E\left[\frac{1}{N}\sum_{k=1}^{N}F(x^k)\right] - \min_{x\in Q}F(x) \le \underbrace{E\left[\frac{1}{N}\sum_{k=1}^{N}F^\mu(x^k)\right] - \min_{x\in Q}F^\mu(x)}_{\operatorname{Regret}_{F^\mu}\left(\{x^k\}_{k=0}^{N-1}\right)\le\varepsilon/2} + \frac{\varepsilon}{2} \le \varepsilon.$$

Таким образом, при условии (19), решение задачи (18) с точностью $\varepsilon/2$ является решением задачи (5) с точностью $\varepsilon$.

## 5. Сведение с помощью техники двойного сглаживания Б.Т. Поляка

Следуя п. 2, будем считать, что $f(x,\xi)$ – выпуклая (вообще говоря, негладкая) функция от $x$, которая удовлетворяет условию Липшица (1).

Возьмём в п. 3
$$F(x,\tilde{\eta}) = f(x+\tau\tilde{e}_1+\mu\tilde{e}_2,\xi), \ F(x) = E_{\tilde{\eta}}\left[F(x,\tilde{\eta})\right], \ \tilde{\eta}=(\xi,\tilde{e}_1,\tilde{e}_2),$$

в качестве $\nabla_x \tilde{F}(x,\eta)$ в МЗС (см. п. 3) возьмём

$$\nabla_x \tilde{F}(x,\eta) = \frac{n}{\mu}\left(\tilde{f}(x+\tau\tilde{e}_1+\mu e_2,\xi) - \tilde{f}(x+\tau\tilde{e}_1,\xi)\right)e_2, \ \eta=(\xi,\tilde{e}_1,e_2), \tag{20}$$

где $\tilde{e}_1 \in RB_2^n(1)$ ($e_2 \in RS_2^n(1)$) – случайный вектор, равномерно распределённый на шаре $B_2^n(1)$ (сфере $S_2^n(1)$), $\tau \le \varepsilon/(4M)$, $\mu \ll \tau$ будут выбраны позже (см. формулу (28), (29)). Будем считать, что $\{\tilde{e}_1^k, e_2^k\}_{k=0}^{N-1}$ независимы в совокупности и независимы от $\{\xi^k\}_{k=0}^{N-1}$.

Обозначим через

$$f^\tau(x) = E_{\tilde{e}_1,\xi}\left[f(x+\tau\tilde{e}_1,\xi)\right], \ f^{\tau,\mu}(x) = E_{\tilde{e}_1,\tilde{e}_2,\xi}\left[f(x+\tau\tilde{e}_1+\mu\tilde{e}_2,\xi)\right] = F(x).$$

**Лемма 1 (см. п. 4).** *Если* $\tau \le \varepsilon/(4M)$, $\mu \le \varepsilon/(4M)$, *то*

$$\operatorname{Regret}_f\left(\{x^k\}_{k=0}^{N-1}\right) \le \operatorname{Regret}_{f^\tau}\left(\{x^k\}_{k=0}^{N-1}\right) + \frac{\varepsilon}{4} \le$$

$$\le \operatorname{Regret}_{f^{\tau,\mu}}\left(\{x^k\}_{k=0}^{N-1}\right) + \frac{\varepsilon}{2} = \operatorname{Regret}_F\left(\{x^k\}_{k=0}^{N-1}\right) + \frac{\varepsilon}{2}. \tag{21}$$

**Лемма 2 (см., например, [3, 5]).** *Пусть*

$$\nabla_x F(x,\eta) = \frac{n}{\mu}\left(f(x+\tau\tilde{e}_1+\mu e_2,\xi) - f(x+\tau\tilde{e}_1,\xi)\right)e_2, \ \eta=(\xi,\tilde{e}_1,e_2),$$

*где* $\tilde{e}_1 \in RB_2^n(1)$, $e_2 \in RS_2^n(1)$ *и* $\tilde{e}_1$, $e_2$, $\xi$ *независимы в совокупности. Тогда*

$$E_\eta\left[\nabla_x F(x,\eta)\right] = \nabla F(x).$$

**Лемма 3 (см. [7]).** *Для последовательности* $\{x^k\}_{k=0}^{N}$, *сгенерированной МЗС (8) с* $\{\nabla_x \tilde{F}(x^k,\eta^k)\}_{k=0}^{N-1}$, *рассчитываемыми по формуле (20), справедлива формула (6) с*

$$\sigma \le \frac{4\delta R\sqrt{n}}{\mu}. \tag{22}$$

В основу доказательства леммы 3 положена лемма 2, оценка (15) и явление концентрации равномерной мере на поверхности евклидовой сферы в малой окрестности экватора [18] (с произвольно выбранным северным полюсом).

**Лемма 4 (И.Н. Усманова, 2015 [19]).** *Пусть* $e_2 \in RS_2^n(1)$, *тогда*

$$E\left[\|e_2\|_q^2\right] \le \min\{q-1, 4\ln n\} \cdot n^{\frac{2}{q}-1} = c_q n^{\frac{2}{q}-1},\ 2 \le q \le \infty, \tag{23}$$

$$E\left[\langle c, e_2\rangle^2 \|e_2\|_q^2\right] \le \frac{4}{3}\|c\|_2^2 \min\{q-1, 4\ln n\} \cdot n^{\frac{2}{q}-2},\ 2 \le q \le \infty. \tag{24}$$

**Лемма 5 (см. [7]).** *Пусть для любого* $e_2 \in S_2^n(1)$

$$\left| f(x+\tau\tilde{e}_1+\mu e_2,\xi) - f(x+\tau\tilde{e}_1,\xi) - \mu\langle\nabla_x f(x+\tau\tilde{e}_1,\xi),e_2\rangle\right| \le \frac{L(\tilde{e}_1,\xi)\mu^2}{2}.$$

*Тогда*

$$E_\eta\left[\|\nabla_x \tilde{F}(x,\eta)\|_q^2\right] \le \frac{3}{4}n^2\mu^2 E_{\tilde{e}_1,\xi}\left[L(\tilde{e}_1,\xi)^2\right] E_{e_2}\left[\|e_2\|_q^2\right] +$$
$$+ 3n^2 E_{\tilde{e}_1,e_2,\xi}\left[\langle\nabla_x f(x+\tau\tilde{e}_1,\xi),e_2\rangle^2 \|e_2\|_q^2\right] + 12\frac{\delta^2 n^2}{\mu^2}E_{e_2}\left[\|e_2\|_q^2\right], \tag{25}$$

*где* $\tilde{e}_1 \in RB_2^n(1)$, $e_2 \in RS_2^n(1)$ *и* $\tilde{e}_1$, $e_2$, $\xi$ *независимы в совокупности*.

Здесь с точностью до почти всюду по $\tilde{e}_1$ под $\nabla_x f(x+\tau\tilde{e}_1,\xi)$ можно понимать обычный (стохастический) градиент, вообще говоря, не гладкой выпуклой функции $f(x,\xi)$, как функции $x$, в точке $x+\tau\tilde{e}_1$. Это следует из теоремы Радемахера [20].

**Лемма 6.** *Пусть* $\tilde{e}_1 \in RB_2^n(1)$ *и* $\tilde{e}_1$, $\xi$ *независимы в совокупности. Тогда*

$$E_{\tilde{e}_1,\xi}\left[L(\tilde{e}_1,\xi)^2\right] \le \frac{16M^2}{3\mu\tau}. \tag{26}$$

Эта оценка следует из того, что все сводится к рассмотрению наиболее неблагоприятного случая в размерности $n=1$ с $f(x,\xi) \equiv M \cdot |x|$ в точке $x = -\mu/2$ с $e_2 = 1$.

Используя (7), (23)–(26) перепишем оценку (25) в виде

$$\tilde{M}^2 = \max_{x \in Q} E_\eta \left[ \left\| \nabla_x \tilde{F}(x,\eta) \right\|_q^2 \right] \leq 4c_q n^{2/q} \left( nM^2 \frac{\mu}{\tau} + M^2 + 3n \frac{\delta^2}{\mu^2} \right). \tag{27}$$

Выберем $\tau$ на пороге (см. лемму 1)

$$\tau = \frac{\varepsilon}{4M}. \tag{28}$$

Выберем $\mu \leq \varepsilon/(4M)$ так, чтобы $nM^2 \cdot \mu/\tau = M^2$, т.е.

$$\mu = \frac{\varepsilon}{4Mn}. \tag{29}$$

Будем считать, что уровень шума таков, что (см. (22), (27))

$$\frac{4\delta R\sqrt{n}}{\mu} \leq \frac{\varepsilon}{4}, \ 3n\frac{\delta^2}{\mu^2} \leq M^2,$$

т.е. условие на максимально допустимый уровень шума $\delta \leq \delta_0$ будет иметь вид

$$\delta_0 \simeq \min\left\{ \frac{\varepsilon^2}{56MRn^{3/2}}, \frac{\varepsilon}{7n^{3/2}} \right\} \simeq \frac{\varepsilon^2}{56MRn^{3/2}}. \tag{30}$$

Интересно сопоставить эту оценку с оценкой из работы [7], в которой дополнительно предполагалась равномерная липшицевость градиента $f(x,\xi)$:

$$\left\| \nabla f(y,\xi) - \nabla f(x,\xi) \right\|_2 \leq L \|y - x\|_2.$$

В таком случае

$$\delta_0 \simeq \min\left\{ \frac{\varepsilon^{3/2}}{16R\sqrt{Ln}}, \frac{M\varepsilon^{1/2}}{\sqrt{96Ln}} \right\} \simeq \frac{\varepsilon^{3/2}}{16R\sqrt{Ln}}.$$

При выбранных таким образом $\tau$, $\mu$ (формулы (28), (29)) и уровне шума $\delta$ из (27), имеем

$$\tilde{M}^2 \leq 12c_q n^{2/q} M^2. \tag{31}$$

Подставляя эту формулу в выражение $N = N(\varepsilon/4)$ (см. формулу (13)), получим число дней (итераций), которые достаточно прожить (сделать), чтобы гарантировать

$$\mathrm{Regret}_F\left(\left\{x^k\right\}_{k=0}^{N-1}\right) \le \varepsilon/2,$$

а, следовательно (см. формулу (21) в лемме 1),

$$\mathrm{Regret}_f\left(\left\{x^k\right\}_{k=0}^{N-1}\right) \le \varepsilon.$$

Резюмируем полученный результат в виде теоремы.

**Теорема.** *Пусть в МЗС из п. 3 (8) (для полного описания метода используется также (11), (31)) подставляется на каждой итерации выражение (20) (сформированное семейной парой из п. 2 на базе доступной информации $\left\{\tilde{f}\left(x_m^k, \xi^k\right), \tilde{f}\left(x_w^k, \xi^k\right)\right\}_k$), где $\tau$, $\mu$ выбраны согласно (28), (29), а уровень шума ограничен $\delta \le \delta_0$ (30). Тогда после*

$$N = \frac{384 c_q n^{2/q} M^2 R^2}{\varepsilon^2}$$

*дней (итераций) имеет место неравенство*

$$\mathrm{Regret}_f\left(\left\{x^k\right\}_{k=0}^{N-1}\right) \le \varepsilon.$$

Другими словами,

$$\mathrm{Regret}_f\left(\left\{x^k\right\}_{k=0}^{N-1}\right) = \tilde{\mathrm{O}}\left(n^{1/q} \frac{MR}{\sqrt{N}}\right).$$

Аналогично устанавливается оценка

$$\mathrm{Regret}_f\left(\left\{x^k\right\}_{k=0}^{N-1}\right) = \mathrm{O}\left(\frac{nM^2 \ln N}{\gamma N}\right)$$

в случае, когда функция $f(x)$ является сильно выпуклой в норме $l_2$ с достаточно большой константой $\gamma > 0$.

Возвращаясь к задаче семейной пары оптимизаторов из п. 2, осталось только заметить, что при условиях (28), (29)

$$\left|\mathrm{Regret}_f\left(\left\{x^k + \tau\tilde{e}_1^k + \mu e_2^k\right\}_{k=0}^{N-1}\right) - \mathrm{Regret}_f\left(\left\{x^k\right\}_{k=0}^{N-1}\right)\right| \le \varepsilon/2,$$

$$\left|\mathrm{Regret}_f\left(\left\{x^k + \tau\tilde{e}_1^k\right\}_{k=0}^{N-1}\right) - \mathrm{Regret}_f\left(\left\{x^k\right\}_{k=0}^{N-1}\right)\right| \le \varepsilon/4.$$

**Литература**